\documentclass[12pt]{article}
\usepackage{amssymb}

\input{tcilatex}
\begin{document}

\begin{center}
\textbf{UNIFORMITY AND SELF-NEGLECTING FUNCTIONS}

\bigskip

\textbf{N. H. BINGHAM and A. J. OSTASZEWSKI}

\bigskip
\end{center}

\noindent \textbf{Abstract. }We relax the continuity assumption in Bloom's
uniform convergence theorem for Beurling slowly varying functions $\varphi $%
. We assume that $\varphi $ has the Darboux property, and obtain results for 
$\varphi $ measurable or having the Baire property. \newline
\noindent \textbf{Keywords}: Karamata slow variation, Beurling slow
variation, Wiener's Tauberian theorem, Beurling's Tauberian theorem, uniform
convergence theorem, Darboux property, Kestelman-Borwein-Ditor theorem,
Baire's category theorem, measurability, Baire property, affine group
action. \newline
\textbf{Classification}:  26A03; 33B99, 39B22, 34D05.\newline

\noindent \textbf{1. Definitions and motivation}\newline

\indent The motivation for this paper may be traced back to two classic
papers. First, in 1930, Karamata [Kar] introduced his theory of \textit{%
regular variation} (in particular, of \textit{slow variation}). Also in 1930
Karamata simplified the Hardy-Littlewood approach to \textit{Tauberian
theorems}; in 1931 he applied his theory of regular variation very
successfully to Tauberian theory. For textbook accounts, see e.g. [BinGT]
(BGT below), Ch. 4, Korevaar [Kor], IV. Secondly, in 1932 Wiener [Wie]
transformed Tauberian theory by working with general kernels (rather than
special kernels as Hardy and Littlewood had done); his method was based on
Fourier transforms. There are two common forms for Wiener's Tauberian
theorem, one for the additive group of reals, one for the multiplicative
group of positive reals. Both concern convolutions: additive convolutions
for the first, with Fourier transforms and Haar measure being Lebesgue
measure $dx$; multiplicative convolutions for the second, with Mellin
transforms and Haar measure $dx/x$. \newline

\noindent \textbf{Theorem W (Wiener's Tauberian theorem)}. \textit{For }$%
K\in L_{1}(\mathbb{R})$\textit{\ with the Fourier transform }$\hat{K}$%
\textit{\ of }$K$\textit{\ non-vanishing on the real line, and }$H\in
L_{\infty }(\mathbb{R})$\textit{: if }%
\[
\int K(x-y)H(y)dy\rightarrow c\int K(y)dy\qquad (x\rightarrow \infty ), 
\]%
\textit{then for all }$G\in L_{1}(\mathbb{R})$\textit{,} 
\[
\int G(x-y)H(y)dy\rightarrow c\int G(y)dy\qquad (x\rightarrow \infty ). 
\]

\indent The corresponding multiplicative form, with $\int_{0}^{\infty
}K(x/y)H(y)dy/y$, is left to the reader. For textbook accounts, see e.g.
Hardy [Har], XII, Widder [Wid], V, BGT Ch. 4, [Kor], II. Usually, in both
regular variation (below) and Tauberian theory, the multiplicative form is
preferred for applications, the additive form for proofs, a practice we
follow here. \newline
\indent The classic summability methods of Ces\`{a}ro and Abel fall easily
into this framework. The next most important family is that of the Euler and
Borel methods; these are tractable by Wiener methods, but are much less
amenable to them; see e.g. [Har] VIII, IX, \S 12.15. Indeed, Tenenbaum
[Ten], motivated by analytic number theory, gives an approach to Tauberian
theory for the Borel method by Hardy-Littlewood rather than Wiener methods. 
\newline
\indent The Borel method is at the root of our motivation here. It is of
great importance, in several areas: analytic continuation by power series
([Har] VIII, IX, Boas [Boa1] \S 5.5); analytic number theory [Ten];
probability theory ([Bin1] -- [Bin4]). \newline
\indent In unpublished lectures, Beurling undertook the task of bringing the
Borel method (and its numerous relatives) within the range of easy
applicability of Wiener methods. His work was later published by Peterson
[Pet] and Moh [Moh]. We state his result as \textit{Beurling's Tauberian
theorem} below, but we must first turn to the two themes of our title. 
\newline
\indent A function $f:\mathbb{R}\rightarrow \mathbb{R}_{+}$ is \textit{%
regularly varying} in Karamata's sense if for some function $g$, 
$$
f(ux)/f(x)\rightarrow g(u)\qquad (x\rightarrow \infty )\qquad \forall \ u>0.%
\eqno(RV)
$$%
It turns out that to get a fruitful theory, one needs \textit{some}
regularity condition on $f$. Karamata himself used continuity. This was
weakened to (Lebesgue) measurability by Korevaar et al. [KvAEdB] in 1949.
Matuszewska [Mat] in 1962 showed that one could also use functions with the
Baire property (briefly: Baire functions). Note that neither of the
measurable and Baire cases contains the other. There are extensive and
useful parallels between the measure and Baire (or category) cases -- see
e.g. Oxtoby [Oxt], [BinO3], [BinO4], [Ost3]; furthermore these have
wide-ranging applications, see e.g. [BinO9], [Ost1], [Ost2] -- particularly
to the Effros Theorem, cf. [Ost4] and \S 5.5. \newline
\indent Subject to a regularity condition ($f$ measurable or Baire, say),
one has: \newline
(i) the \textit{uniform convergence theorem} (UCT): $(RV)$ holds \textit{%
uniformly} on compact $u$-sets; \newline
(ii) the \textit{characterization theorem}: $g(u)=u^{\rho }$ for some $\rho $%
, called the \textit{index} of regular variation. \newline
The class of such $f$, those regularly varying with index $\rho $, is
written $R_{\rho }$. One can reduce to the case $\rho =0$, of the (Karamata) 
\textit{slowly varying} functions, $R_{0}$. Then 
$$
f(xu)/f(x)\rightarrow 1\qquad (x\rightarrow \infty )\qquad \forall \ u\eqno%
(SV)
$$%
working multiplicatively, or 
$$
h(x+u)-h(x)\rightarrow 0\qquad (x\rightarrow \infty )\qquad \forall \ u\eqno%
(SV_{+})
$$%
working additively; either way the convergence is uniform on compact $u$%
-sets (in the line for $(SV)$, the half-line for $(SV_{+})$). \newline
\indent Beurling observed in his lectures that the function $\sqrt{x}$ --
known to be crucial to the Tauberian theory of the Borel method -- has a
property akin to Karamata's slow variation. We say that $\varphi >0$ is 
\textit{Beurling slowly varying}, $\varphi \in BSV$, if $\varphi (x)=o(x)$
as $x\rightarrow \infty $ and 
$$
\varphi (x+t\varphi (x))/\varphi (x)\rightarrow 1\qquad (x\rightarrow \infty
)\qquad \forall \ t.\eqno(BSV)
$$%
Using the additive notation $h:=\log \varphi $ (whenever convenient):%
$$
h(x+t\varphi (x))-h(x)\rightarrow 0\qquad (x\rightarrow \infty )\qquad
\forall \ t.\eqno(BSV_{+})
$$%
If (as in the UCT for Karamata slow variation) the convergence here is
locally uniform in $t$, we say that $\varphi $ is \textit{self-neglecting}, $%
\ \varphi \in SN$; we write $(SN)$ for the corresponding strengthening of $%
(BSV)$. \newline
\indent We may now state Beurling's extension to Wiener's Tauberian theorem
(for background and further results, see \S 5.1 below and [Kor]). \newline

\noindent \textbf{Theorem (Beurling's Tauberian theorem)}.\textit{\ If }$%
\varphi \in BSV$\textit{, }$K\in L_{1}(\mathbb{R})$\textit{\ with }$\hat{K}$%
\textit{\ non-zero on the real line, }$H$\textit{\ is bounded, and }%
\[
\int K\Bigl(\frac{x-y}{\varphi (x)}\Bigr)H(y)dy/\varphi (x)\rightarrow c\int
K(y)dy\qquad (x\rightarrow \infty ), 
\]%
\textit{then for all }$G\in L_{1}(\mathbb{R})$\textit{,} 
\[
\int G\Bigl(\frac{x-y}{\varphi (x)}\Bigr)H(y)dy/\varphi (x)\rightarrow c\int
G(y)dy\qquad (x\rightarrow \infty ). 
\]

\indent Notice that the arguments of $K$ and $G$ here involve \textit{both}
the additive group operation on the line \textit{and} the multiplicative
group operation on the half-line. Thus Beurling's Tauberian theorem,
although closely related to Wiener's (which it contains, as the case $%
\varphi \equiv 1$), is structurally different from it. One may also see here
the relevance of the \textit{affine group}, already well used for regular
variation (see e.g. BGT \S 8.5.1, [BinO2], and \S 3 below). \newline
\indent Analogously to Karamata's UCT, the following result was proved by
Bloom in 1976 [Blo]. A slightly extended and simplified version is in BGT,
Th. 2.11.1. \newline

\noindent \textbf{Theorem (Bloom's theorem)}. \textit{If }$\varphi \in BSN$%
\textit{\ with }$\varphi $\textit{\ continuous, then }$\varphi \in SN$%
\textit{: }$(BSN)$\textit{\ holds locally uniformly.}\newline

\indent The question as to whether one can extend this to $\varphi $
measurable and/or Baire has been open ever since; see BGT \S 2.11, [Kor]
IV.11 for textbook accounts. Our purpose here is to give some results in
this direction. This paper is part of a series (by both authors, and by the
second author, alone and in [MilO] with Harry I. Miller) on our new theory
of \textit{topological regular variation}; see e.g. [BinO1-11], [Ost1-4] and
the references cited there. One of the objects achieved was to find the
common generalization of the measurable and Baire cases. This involves
infinite combinatorics, in particular such results as the
Kestelman-Borwein-Ditor theorem (KBD -- see e.g. [MilO]), the category
embedding theorem [BinO4] (quoted in \S 3 below) and shift-compactness
[Ost3]. A by-product was the realization that, although the Baire case came
much later than the measurable case, it is in fact the more important. One
can often handle both cases together bitopologically, using the Euclidean
topology for the Baire case and the density topology for the measurable
case; see \S 3 below and also [BinO4]. Such measure-category duality only
applies to \textit{qualitative} measure theory (where all that counts is
whether the measure of a set is zero or positive, not its numerical value).
We thus seek to avoid \textit{quantitative} measure-theoretic arguments; see 
\S 5.5.\newline
\indent Our methods of proof (as with our previous studies in this area)
involve tools from infinite combinatorics, and replacement of quantitative
measure theory by qualitative measure theory.

\bigskip

\noindent \textbf{2. Extensions of Bloom's Theorem: Monotone functions}%
\newline

\indent We suggest that the reader cast his eye over the proof of Bloom's
theorem, in either [Blo] or BGT \S 2.11 -- it is quite short. Like most
proofs of the UCT for Karamata slow variation, it proceeds by contradiction,
assuming that the desired uniformity fails, and working with two sequences, $%
t_{n}\in \lbrack -T,T]$ and $x_{n}\rightarrow \infty $, witnessing to its
failure. \newline
\indent The next result, in which we assume $\varphi $ monotone ($\varphi $
increasing to infinity is the only case that requires proof) is quite
simple. But it is worth stating explicitly, for three reasons: \newline
1. It is a complement to Bloom's theorem, and to the best of our knowledge
the first new result in the area since 1976. \newline
2. The case $\varphi $ increasing is by far the most important one for
applications. For, taking $G$ the indicator function of an interval in
Beurling's Tauberian theorem, the conclusion there has the form of a \textit{%
moving average}: 
\[
\frac{1}{a\varphi (x)}\int_{x}^{x+a\varphi (x)}H(y)dy\rightarrow c\qquad
(x\rightarrow \infty )\qquad \forall a\ >0.
\]%
Such moving averages are \textit{Riesz (typical) means} and here $\varphi $
increasing to $\infty $ is natural in context. For a textbook account, see
[ChaM]; for applications, in analysis and probability theory, see [Bin5],
[BinG1], [BinG2], [BinT]. The prototypical case is $\varphi (x)=x^{\alpha }$ 
$(0<\alpha <1)$; this corresponds to $X\in L_{1/\alpha }$ for the
probability law of $X$.\newline
3. Theorem 1 below is closely akin to results of de Haan on the Gumbel law $%
\Lambda $ in extreme-value theory; see \S 5.8 below.\newline
\indent We offer three proofs (two here and a third after Theorem 2M) of the
result, as each is short and illuminating in its own way.\newline
\indent For the first, recall that if a sequence of monotone functions
converges pointwise to a continuous limit, the convergence is uniform on
compact sets. See e.g. P\'{o}lya and Szeg\H{o} [PolS], Vol. 1, p.63, 225,
Problems II 126, 127, Boas [Boa2], \S 17, p.104-5. (The proof is a simple
compactness argument. The result is a complement to the better-known result
of Dini, in which it is the convergence, rather than the functions, that is
monotone; see e.g. [Rud], 7.13.)\newline

\noindent \textbf{Theorem 1 (Monotone Beurling UCT)}. \textit{If }$\varphi
\in BSV$\textit{\ is monotone, }$\varphi \in SN$\textit{: the convergence in 
}$(BSV)$\textit{\ is locally uniform. }\newline
\bigskip

\noindent \textit{First proof}. As in [Blo] or BGT \S 2.11, we proceed by
contradiction. Pick $T>0$, and assume the convergence is not uniform on $%
[-T,0]$ (the case $[0,T]$ is similar). Then there exists $\varepsilon
_{0}\in (0,1)$, $t_{n}\in \lbrack -T,0]$ and $x_{n}\rightarrow \infty $ such
that 
\[
|\varphi (x_{n}+t_{n}\varphi (x_{n}))/\varphi (x_{n})-1|\geq \varepsilon
_{0}\qquad \forall \ n. 
\]%
Write 
\[
f_{n}(t):=\varphi (x_{n}+t\varphi (x_{n}))/\varphi (x_{n})-1. 
\]%
Then $f_{n}$ is monotone, and tends pointwise to 0 by $(BSV)$. So by the P%
\'{o}lya-Szeg\H{o} result above, the convergence is uniform on compact sets.
This contradicts $|f_{n}(t_{n})|\geq \varepsilon _{0}$ for all $n$. $\square 
$

\bigskip

\indent The second proof is based on the following result, thematic for the
approach followed in \S 4. We need some notation that will also be of use
later. Below, \ $x>0$ will be a continuous variable, or a sequence $%
x:=\{x_{n}\}$ diverging to $+\infty $ (briefly, \textit{divergent sequence}%
), according to context. We put%
\[
V_{n}^{x}(\varepsilon ):=\{t\geq 0:|\varphi (x_{n}+t\varphi (x_{n}))/\varphi
(x_{n})-1|\leq \varepsilon \},\text{ }H_{k}(\varepsilon
):=\bigcap\nolimits_{n\geq k}V_{n}^{x}(\varepsilon ). 
\]%
\newline
\noindent \textbf{Lemma 1. }\textit{For }$\varphi >0$\textit{\ monotonic
increasing and }$\{x_{n}\}$ \textit{a divergent sequence,\ each set }$%
V_{n}^{x}(\varepsilon ),$ \textit{and so also each set }$H_{k}^{x}(%
\varepsilon ),$ \textit{is an interval containing }$0$.\newline

\bigskip

\noindent \textit{Proof of Lemma 1}. Since $x+s\varphi (x)>x$ for $s>0,$ one
has $1\leq \varphi (x+t\varphi (x))/\varphi (x).$ Also if $0<s<t,$ then, as $%
\varphi (x)>0$, one has $x+s\varphi (x)<x+t\varphi (x).$ So if $t\in
V_{n}^{x}(\varepsilon ),$ then%
\[
1\leq \varphi (x_{n}+s\varphi (x_{n}))/\varphi (x_{n})\leq \varphi
(x_{n}+t\varphi (x_{n}))/\varphi (x_{n})\leq 1+\varepsilon ,
\]%
and so $s\in V_{n}^{x}(\varepsilon ).$ The remaining assertions now follow,
because an intersection of intervals containing $0$ is an interval
containing $0$. $\square $

\bigskip

\noindent \textit{Second proof of Theorem 1}. Suppose otherwise; then there
are $\varepsilon _{0}>0$ and sequences $x_{n}:=x(n)\rightarrow \infty $ and $%
u_{n}\rightarrow u_{0}$ such that%
\begin{equation}
|\varphi (x_{n}+u_{n}\varphi (x_{n}))/\varphi (x_{n})-1|\geq \varepsilon
_{0},\text{ }(\forall n\in \mathbb{N}).  \tag{all}  \label{all}
\end{equation}%
Since $\varphi $ is Beurling slowly varying the increasing sets $%
H_{k}^{x}(\varepsilon _{0})$ cover $\mathbb{R}_{+}$ and so, being increasing
intervals (by Lemma 1), their interiors cover the compact set $%
K:=\{u_{n}:n=0,1,2,...\}.$ So for some integer $k$ the set $%
H_{k}^{x}(\varepsilon _{0})$ already covers $K,$ and then so does $%
V_{k}^{x}(\varepsilon _{0}).$ But this implies that 
\[
|\varphi (x_{k}+u_{k}\varphi (x_{k}))/\varphi (x_{k})-1|<\varepsilon _{0},
\]%
contradicting (all) at $n=k.$ $\square $

\bigskip

\noindent \textit{Remark}. Of course the uniformity property of $\varphi $
is equivalent to the sets $H_{k}^{x}(\varepsilon )$ containing arbitrarily
large intervals $[0,t]$ for large enough $k$ (for all divergent $\{x_{n}\}).$

\bigskip

\noindent \textbf{3. Combinatorial preliminaries}\newline

\indent We work in the affine group $\mathcal{A}ff$ acting on $(\mathbb{R},+)
$ using the notation 
\[
\gamma _{n}(t)=c_{n}t+z_{n},
\]%
where $c_{n}\rightarrow c_{0}=c>0$ and $z_{n}\rightarrow 0$ as $n\rightarrow
\infty $. These are to be viewed as (self-) homeomorphisms of $\mathbb{R}$
under either $\mathcal{E},$ the Euclidean topology, or $\mathcal{D}$, the
Density topology. Recall that the open sets of $\mathcal{D}$ are measurable
subsets, all points of which are (Lebesgue) density points, and that (i)
Baire sets under $\mathcal{D}$ are precisely the Lebesgue measurable sets,
(ii) the nowhere dense sets of $\mathcal{D}$ are precisely the null sets,
and (iii) Baire's Theorem holds for $\mathcal{D}$. (See Kechris [Kech]
17.47.) Below we call a set \textit{negligible} if, according to the
topological context of $\mathcal{E}$ or $\mathcal{D}$, it is meagre/null. A
property holds for `quasi all' elements of a set if it holds for all but a
negligible subset. We recall the following definition and Theorem from
[BinO4], which we apply taking the space $X$ to be $\mathbb{R}$ with one of $%
\mathcal{E}$ or $\mathcal{D}$.\newline
\noindent \textbf{Definition. }A sequence of homeomorphisms $%
h_{n}:X\rightarrow X$ satisfies the \textit{weak category convergence }%
condition (wcc) if:\newline
\indent For any non-meagre open set $U\subseteq X,$ there is a non-meagre
open set $V\subseteq U$ such that for each $k\in N,$%
\[
\bigcap\nolimits_{n\geq k}V\backslash h_{n}^{-1}(V)\text{ is meagre.}
\]

\bigskip

\noindent \textbf{Theorem CET (Category Embedding Theorem).}\textit{\ Let }$%
X $\textit{\ be a topological space and }$h_{n}:X\rightarrow X$\textit{\ be
homeomorphisms satisfying (wcc). Then for any Baire set }$T,$\textit{\ for
quasi-all }$t\in T$\textit{\ there is an infinite set} $\mathbb{M}%
_{t}\subseteq \mathbb{N}$ \textit{such that}%
\[
\{h_{m}(t):m\in \mathbb{M}_{t}\}\subseteq T. 
\]%
\newline
From here we deduce: \bigskip

\noindent \textbf{Lemma 2 (Affine Two-sets Lemma). }\textit{For }$%
c_{n}\rightarrow c>0$ \textit{and }$z_{n}\rightarrow 0,$ \textit{if }$%
cB\subseteq A$ \textit{for }$A,B$ \textit{non-negligible (measurable/Baire),
then for quasi all }$b\in B$ \textit{there exists an infinite set }$\mathbb{M%
}=\mathbb{M}_{b}\subseteq \mathbb{N}$ \textit{such that}%
\[
\{\gamma _{m}(b)=c_{m}b+z_{m}:m\in \mathbb{M}\}\subseteq A.
\]

\noindent \textit{Proof}.\textbf{\ }It is enough to prove the existence of
one such point $b,$ as the Generic Dichotomy Principle (for which see
[BinO7, Th. 3.3]) applies here, because we may prove existence of such a $b$
in any non-negligible $\mathcal{G}_{\delta }$-subset $B^{\prime }$ of $B$,
by replacing $B$ below with $B^{\prime }$. (One checks that the set of $b$s
with the desired property is Baire, and so its complement in $B$ cannot
contain a non-negligible $\mathcal{G}_{\delta }$.)\newline
\indent Writing $T:=cB$ and $w_{n}=c_{n}c^{-1},$ so that $c_{n}=w_{n}c$ and $%
w_{n}\rightarrow 1,$ put%
\[
h_{n}(t):=w_{n}t+z_{n}.
\]%
Then $h_{n}$ converges to the identity in the supremum metric, so (wcc)
holds by Th. 6.2 of [BinO6] (First Verification Theorem) and so Theorem CET
above applies for the Euclidean case; applicability in the measure case is
established as Cor. 4.1 of [BinO2]. (This is the basis on which the affine
group preserves negligibility.) So there are $t\in T$ and an infinite set of
integers $\mathbb{M}$ with%
\[
\{w_{m}t+z_{m}:m\in \mathbb{M}\}\subseteq T.
\]%
But $t=cb$ for some $b\in B\ $and so, as $w_{m}c=c_{m},$ one has%
\[
\{c_{m}b+z_{m}:m\in \mathbb{M}\}\subseteq cB\subseteq A.\qquad \square 
\]

\bigskip

\noindent \textbf{4. Extensions to Bloom's theorem: Darboux property} 
\newline

\indent In this section we generalize Bloom's Theorem and simplify his
proof. Bloom uses continuity only through the intermediate value property --
that if a (real-valued) function attains two values, it must attain all
intermediate values. This is the \textit{Darboux property}. It is much
weaker than continuity -- it does not imply measurability, nor the Baire
property. For measurability, see the papers of Halperin [Halp1,2]; for the
Baire property, see e.g. [PorWBW] and also \S 5.3 below.\newline
\indent We use Lemma 2 above to prove Theorem 3B below, which implies
Bloom's Theorem, as continuous functions are Baire and have the Darboux
property. We note a result of Kuratowski and Sierpi\'{n}ski [KurS] that for
a function of Baire class 1 (for which see \S 5.3) the Darboux property is
equivalent to its graph being connected; so Theorem 2 goes beyond the class
of functions considered by Bloom.

We begin with some infinite combinatorics associated with a positive
function $\varphi \in BSV$.

\bigskip

\textbf{Definitions. }Say that $\{u_{n}\}$ with limit $u$ is a \textit{%
witness sequence at }$u$\textit{\ }(for non-uniformity in $\varphi $) if
there are $\varepsilon _{0}>0$ and a divergent sequence $x_{n}$ such that
for $h=\log \varphi $ 
\begin{equation}
|h(x_{n}+u_{n}\varphi (x_{n}))-h(x_{n})|>\varepsilon _{0}\qquad \forall \
n\in \mathbb{N}.  \tag{$\varepsilon _{0}$}
\end{equation}%
Say that $\{u_{n}\}$ with limit $u$ is a \textit{divergent witness sequence}
if also 
\[
h(x_{n}+u_{n}\varphi (x_{n}))-h(x_{n})\rightarrow \pm \infty .
\]%
\textit{\ }Thus a divergent witness sequence is a special type of witness
sequence, but, as we show, it is these that characterize absence of
uniformity in the class $BSV$. 

We begin with a lemma that yields simplifications later; it implies a
Beurling analogue of the Bounded Equivalence Principle in the Karamata
theory, first noted in [BinO1]. As it shifts attention to the origin, we
call it the Shift Lemma. Below \textit{uniform near} a point $u$ means
`uniformly on sequences converging to $u$' and is equivalent to local
uniformity at $u$ (i.e. on compact neighbourhoods of $u).$

\bigskip 

\noindent \textbf{Lemma 3 (Shift Lemma).}\textit{\ For any }$u,$\textit{\
convergence in }$(BSV_{+})$\textit{\ is uniform near }$t=0$\textit{\ iff it
is uniform near }$t=u.$

\bigskip

\noindent \textit{Proof.} Take $z_{n}\rightarrow 0.$ For any $u$ write $%
y_{n}:=x_{n}+u\varphi (x_{n}),$ $\gamma _{n}:=\varphi (x_{n})/\varphi (y_{n})
$ and $w_{n}=\gamma _{n}z_{n}$. Then $\gamma _{n}\rightarrow 1,$ so $%
w_{n}\rightarrow 0.$ But%
\begin{eqnarray*}
&&h(x_{n}+(u+z_{n})\varphi (x_{n}))-h(x_{n}) \\
&=&[h(x_{n}+u\varphi (x_{n})+z_{n}\varphi (x_{n}))-h(x_{n}+u\varphi
(x_{n}))]+[h(x_{n}+u\varphi (x_{n}))-h(x_{n})] \\
&=&[h(y_{n}+z_{n}\gamma _{n}\varphi (y_{n}))-h(y_{n})]+[h(y_{n})-h(x_{n})],
\end{eqnarray*}
i.e.%
\[
h(x_{n}+(u+z_{n})\varphi (x_{n}))-h(x_{n})=[h(y_{n}+w_{n}\varphi
(y_{n}))-h(y_{n})]+[h(y_{n})-h(x_{n})].
\]%
The result follows, since $h(y_{n})-h(x_{n})\rightarrow 0,$ as $\varphi \in
BSV$ and $y_{n}\rightarrow \infty $. $\square $

\bigskip

\noindent \textbf{Theorem 2B (Divergence Theorem -- Baire version).} \textit{%
If }$\varphi \in BSV$\textit{\ has the Baire property and }$u_{n}$\textit{\
with limit }$u$ \textit{is a witness sequence, then }$u_{n}$\textit{\ is a
divergent witness sequence.}

\bigskip

\noindent \textit{Proof.} As $u_{n}$ is a witness sequence, for some $%
x_{n}\rightarrow \infty $ and $\varepsilon _{0}>0$ one has ($\varepsilon _{0}
$), with $h=\log \varphi ,$ as always. By the Shift Lemma (Lemma 3) we may
assume that $u=0.$ So (as in the Proof of Lemma 2) we will write $z_{n}$ for 
$u_{n}.$ If $z_{n}$ is not a divergent witness sequence, then $\{{\varphi
(x_{n}+}${$z$}$_{n}{\varphi (x_{n}))}/{\varphi (x_{n})\}}$ contains a
bounded subsequence and so a convergent sequence. W.l.o.g. we thus also have%
\begin{equation}
c_{n}:={\varphi (x_{n}+z}_{n}{\varphi (x_{n}))}/{\varphi
(x_{n})\longrightarrow c\in (0,\infty ).}  \tag{lim}
\end{equation}%
Write $\gamma _{n}(s):=c_{n}s+z_{n}$ and $y_{n}:=x_{n}+z_{n}\varphi (x_{n}).$
Then $y_{n}=x_{n}(1+z_{n}\varphi (x_{n})/x_{n})\rightarrow \infty $ and 
\begin{equation}
\left\vert h{(y}_{n}{)-h(x_{n})}\right\vert \geq \varepsilon _{0}. 
\tag{$\varepsilon _{0}^{\prime }$}
\end{equation}%
\indent Now take $\eta =\varepsilon _{0}/3$ and amend the notation of \S 2
to read%
\[
V_{n}^{x}(\eta ):=\{s\geq 0:|h(x_{n}+s\varphi (x_{n}))-h(x_{n})|\leq \eta \},%
\text{ }H_{k}^{x}(\eta ):=\bigcap\nolimits_{n\geq k}V_{n}^{x}(\eta ).
\]
These are Baire sets, and 
$$
\mathbb{R}=\bigcup\nolimits_{k}H_{k}^{x}(\eta
)=\bigcup\nolimits_{k}H_{k}^{y}(\eta ),\eqno\text{(cov)}
$$%
as $\varphi \in BSV$. The increasing sequence of sets $\{H_{k}^{x}(\eta )\}$
covers $\mathbb{R}.$ So for some $k$ the set $H_{k}^{x}(\eta )$ is
non-negligible. Furthermore, as $c>0,$ the set $c^{-1}H_{k}^{x}(\eta )$ is
non-negligible and so, by (cov), for some $l$ the set%
\[
B:=(c^{-1}H_{k}^{x}(\eta ))\cap H_{l}^{y}(\eta )
\]%
is also non-negligible. Taking $A:=H_{k}^{x}(\eta ),$ one has $B\subseteq
H_{l}^{y}(\eta )$ and $cB\subseteq A$ with $A,B\ $non-negligible. Applying
Lemma 2 to the maps $\gamma _{n}(s)=c_{n}s+z_{n}$, there exist $b\in B$ and
an infinite set\textit{\ }$\mathbb{M}$ such that%
\[
\{c_{m}b+z_{m}:m\in \mathbb{M}\}\subseteq A=H_{k}^{x}(\eta ).
\]%
That is, as $B\subseteq H_{l}^{y}(\eta ),$ there exist $t\in H_{l}^{y}(\eta )
$ and an infinite $\mathbb{M}_{t}$ such that%
\[
\{\gamma _{m}(t)=c_{m}t+z_{m}:m\in \mathbb{M}_{t}\}\subseteq H_{k}^{x}(\eta
).
\]%
In particular, for this $t$ and $m\in \mathbb{M}_{t}$ with $m>k,l$ one has  
\[
t\in V_{m}^{y}(\eta )\hbox{ and }\gamma _{m}(t)\in V_{m}^{x}(\eta ).
\]%
Fix such an $m.$ As $\gamma _{m}(t)\in V_{m}^{x}(\eta ),$%
\begin{equation}
\left\vert h{(x_{m}+\gamma }_{m}({t)\varphi (x_{m}))}-h{(x_{m})}\right\vert
\leq \eta .  \tag{*}  \label{*}
\end{equation}%
But $\gamma _{m}(t)=c_{m}t+z_{m}=z_{m}+t\varphi (y_{m})/\varphi (x_{m}),$ so 
\[
x_{m}+\gamma _{m}(t)\varphi (x_{m})=x_{m}+z_{m}\varphi (x_{m})+t\varphi
(y_{m})=y_{m}+t\varphi (y_{m}),
\]%
`absorbing' the affine shift ${\gamma }_{m}(${$t$}${)}$ into $y.$ So, by
(*), 
\[
\left\vert h{(y_{m}+t\varphi (y_{m}))}-h{(x_{m})}\right\vert \leq \eta .
\]%
But $t\in V_{m}^{y}(\eta ),$ so%
\[
\left\vert h{(y_{m}+t\varphi (y_{m}))}-h{(y_{m})}\right\vert \leq \eta .
\]%
Using the triangle inequality and combining the last two inequalities, we
have%
\begin{eqnarray*}
|h(y_{m})-h(x_{m})| &\leq &|h(y_{m}+t\varphi
(y_{m}))-h(y_{m})|+|h(y_{m}+t\varphi (y_{m}))-h(x_{m})| \\
&\leq &2\eta <\varepsilon _{0},
\end{eqnarray*}%
contradicting ($\varepsilon _{0}^{\prime }$). \U{25a1}\newline

\bigskip

\noindent \textbf{Theorem 2M (Divergence Theorem -- Measure version).} 
\textit{If }$\varphi \in BSV$\textit{\ is measurable and }$u_{n}$\textit{\
with limit }$u$ \textit{is a witness sequence, then }$u_{n}$\textit{\ is a
divergent witness sequence.}

\bigskip

\noindent \textit{Proof}. The argument above applies, with the density
topology $\mathcal{D}$ in place of the Euclidean topology $\mathcal{E}$ (the
real line is still a Baire space, as remarked earlier). \U{25a1}\newline

\bigskip

As an immediate corollary we have:

\bigskip

\noindent \textit{Third Proof of Theorem 1}. If not, then there exists a
witness sequence $u_{n}$ with limit $u.$ By Lemma 3, w.l.o.g. $u>0.$ Let $%
v>u>w>0.$ Since $\varphi \in BSV,$ 
\[
\varphi (x_{n}+v\varphi (x_{n}))/\varphi (x_{n})\rightarrow 1\text{ and }%
\varphi (x_{n}+w\varphi (x_{n}))/\varphi (x_{n})\rightarrow 1,
\]%
so there is $N$ such that both $(1/2)\varphi (x_{n})<\varphi (x_{n}+w\varphi
(x_{n}))$ and  $\varphi (x_{n}+v\varphi (x_{n}))<2\varphi (x_{n})$ for all $%
n>N.$ By increasing $N$ if necessary we may assume that $w<u_{n}<v$ for $n>N.
$ But then%
\[
(1/2)\varphi (x_{n})<\varphi (x_{n}+w\varphi (x_{n}))<\varphi
(x_{n}+u_{n}\varphi (x_{n}))<\varphi (x_{n}+v\varphi (x_{n}))<2\varphi
(x_{n})
\]%
implies that $1/2<\varphi (x_{n}+u_{n}\varphi (x_{n}))/\varphi (x_{n})<2,$
contradicting Theorem 2B/2M, as $\varphi $ is Baire/measurable. \U{25a1}

\bigskip

We now deduce

\bigskip

\noindent \textbf{Theorem 3B (Beurling-Darboux UCT: Baire version).} \textit{%
If }$\varphi \in BSV$\textit{\ has the Baire and Darboux properties, then }$%
\varphi \in SN$\textit{: }$(BSV)$\textit{\ holds locally uniformly. }\newline

\bigskip

\noindent \textit{Proof.} Suppose the conclusion of the theorem is false.
Take $h=\log \varphi .$ Then there exists a witness sequence $v_{n}$ with
limit $v$ and in particular for some $x_{n}\rightarrow \infty $ and $%
\varepsilon _{0}>0$ one has the inequality ($\varepsilon _{0})$ above
modified so that $v_{n}$ replaces $u_{n}.$ 

We construct below a convergent sequence $u_{n}$, with limit $u$ say, such
that%
\begin{equation}
c_{n}:=h{(x_{n}+u}_{n}{\varphi (x_{n}))}-h{(x_{n})\longrightarrow c\in
(-\infty ,\infty ),}  \tag{lim+}
\end{equation}%
and also the unmodified ($\varepsilon _{0}$) holds. This will contradict
Theorem 2B.\newline
The proof here splits according as {$h$}${(x_{n}+v}_{n}{\varphi (x_{n}))}-${$%
h$}${(x_{n})}$ are bounded.

Case (i) The differences {$h$}${(x_{n}+v}_{n}{\varphi (x_{n}))}-${$h$}${%
(x_{n})}$ diverge to $\pm \infty $. Here we appeal to the Darboux property
to replace the sequence $\{v_{n}\}$ with another sequence $\{u_{n}\}$ for
which the corresponding differences are convergent.\newline
\indent Now $f_{n}(t)=h(x_{n}+t\varphi (x_{n}))-h(x_{n})$ has the Darboux
property and $f_{n}(0)=0.$ Either $f_{n}(v_{n})\geq \varepsilon _{0}$ and so
there exists $u_{n}$ between $0$ and $v_{n}$ with $f_{n}(u_{n})=\varepsilon
_{0}$, or $-f_{n}(v_{n})\geq \varepsilon _{0}$, and so there exists $u_{n}$
with $-f_{n}(u_{n})=\varepsilon _{0}$. Either way $|f_{n}(u_{n})|=%
\varepsilon _{0}$. W.l.o.g. $\{u_{n}\}$ is convergent with limit $u$ say,
since $\{v_{n}\}$ is so, and now (lim+) and ($\varepsilon _{0}$) hold, the
latter as in fact%
\[
\left\vert {h(x_{n}+u}_{n}{\varphi (x_{n}))}-{h(x_{n})}\right\vert
=\varepsilon _{0}.
\]

Case (ii) {$h$}${(x_{n}+v}_{n}{\varphi (x_{n}))}-${$h$}${(x_{n})}$ are
bounded. In this case we can get (lim+) by passing to a subsequence.\newline
\indent In either case we contradict Theorem 2B. \U{25a1}

\bigskip

\noindent \textbf{Theorem 3M (Beurling-Darboux UCT - Measure version).} 
\textit{If }$\varphi \in BSV$\textit{\ is measurable and has the Darboux
property, then }$\varphi \in SN$\textit{: }$(BSV)$\textit{\ holds locally
uniformly.} \newline

\noindent \textit{Proof}. The argument above applies, appealing this time to
Theorem 2M. \U{25a1}\newline

\noindent \textit{Remarks}. 1. The Darboux property in Theorems 3 above may
be replaced with a weaker local property. It is enough to require that $%
\varphi $ be \textit{locally range-dense} -- i.e. that at each point $t$
there is a bounded open neighbourhood $I_{t}$ such that the range $\varphi
\lbrack I_{t}]$ is dense in the interval $(\inf \varphi \lbrack I_{t}],\sup
\varphi \lbrack I_{t}])$ -- or be in the class $\mathfrak{A}_{0}$ of [BruC, 
\S 2], cf. also [BruCW].

\noindent 2. The proofs of Theorems 3B and 3M begin as Bloom's does, but
only in the case (i) of the first step, and even then we appeal to the
Darboux property rather than the much stronger assumption of continuity.
Thereafter, we are able to use Theorem CET to base the rest of the proof on
Baire's category theorem. This enables us to handle Theorems 3B and 3M
together, by qualitative measure theory; see the end of \S 1 and \S 5.5
below. By contrast, the proofs of Bloom's theorem in [Blo] and BGT \S 2.11
use quantitative measure theory; see \S 5.5.\newline

\bigskip

\noindent \textbf{5. Complements}\newline
\noindent 1. \textit{Beurling's Tauberian theorem: approximation form. }%
Recall (see e.g. [Kor] II.8) that Wiener's Tauberian theorem is a
consequence of Wiener's approximation theorem: that for $f\in L_{1}(\mathbb{R%
})$ the following are equivalent:

\noindent (i) linear combinations of translates of $f$ are dense in $L_{1}(%
\mathbb{R}),$\newline
\noindent (ii) the Fourier transform $\hat{f}$ of $f$ has no real zeros.%
\newline
\indent The result is the key to Beurling's Tauberian theorem ([Kor] IV Th.
11.1). Rate of convergence results (Tauberian remainder theorems) are also
possible; see e.g. [FeiS], [Kor] VII.13).\newline
\indent The theory extends to Banach algebras (indeed, played a major role
in their development). In this connection ([Kor], V.4) we mention weighted
versions of $L_{1}:$ for \textit{Beurling weights} $\omega $ -- positive
measurable functions on $\mathbb{R}$ with subadditive logarithms, $\omega
(t+u)\leq \omega (t)\omega (u)$ -- define $L_{\omega }=L_{1,\omega }$ to be
the set of $f$ with%
\[
||f||=||f||_{1,\omega }:=\int_{\mathbb{R}}|f(t)|\omega (t)dt<\infty .
\]%
Then (\textit{Beurling's approximation theorem}): Wiener's approximation
theorem extends to the weighted case when $\omega $ satisfies the
nonquasi-analyticity condition%
\[
\int_{\mathbb{R}}\frac{|\log \omega (t)|}{1+t^{2}}dt<\infty .
\]%
This condition has been extensively studied (see e.g. [Koo]) and is
important in probability theory (work by Szeg\H{o} -- see e.g. [Bin6]).%
\newline
\noindent 2. \textit{Representation. }As with Karamata slow variation,
Beurling slow variation has a representation theorem: $\varphi \in SN$ iff $%
\varphi >0$ and%
\[
\varphi (x)=c(x)\int_{0}^{x}e(u)du,
\]%
with $e(.)\rightarrow 0,c(.)\rightarrow c\in (0,\infty );$ as in BGT \S %
2.11, [BinO5, Part II] we may take $e(.)\in \mathcal{C}^{\infty }$ (so the
integral is smooth), and then $c(.)$ has the same degree of regularity
(Baire/measurable, descriptive character, etc.) as $\varphi (.).$ The
treatment of BGT \S 2.11 goes over to the setting here without change. So
too does the drawback that the representation on the right above does not
necessarily imply that $\varphi $ is positive -- this has to be assumed, or
to be given from context.\newline
\noindent 3. \textit{Functions of Baire class 1. }Recall that Baire class 1
functions -- briefly, Baire-1 functions -- are limits of sequences of
continuous functions, and then the Baire hierarchy is defined by successive
passages to the limit. See e.g. Bruckner and Leonard [BruL] \S 2, and the
extensive bibliography given there, [Sol]. Compare [Kech] \S 24.B. The union
of the classes in the Baire hierarchy gives the Borel functions; see e.g.
[Nat] Ch. XV. Since Borel functions are (Lebesgue) measurable and Baire
(have the Baire property)\footnote{%
In general one needs to distinguish between Borel and Baire measurability
(cf. [Halm, \S 51] and [BinO5, \S 11]), \ but the two coincide for real
analysis, our context here -- see [Kech, 24.3].}, the Baire-1 functions are
both measurable and Baire (see \ e.g. [Kur] \S 11).\newline
\indent Lee, Tang and Zhao [LeeTZ] define a concept of \textit{weak
separation} involving neighbourhood assignments. They show that for
real-valued functions on a Polish space this is equivalent to being of Baire
class 1. Their result is greatly generalized by Bouziad [Bou].\newline
\noindent 4. \textit{Darboux functions of Baire class 1. }We recall that a
Darboux function need be neither measurable nor (with the property of) Baire
-- hence the need to impose Darboux-Lebesgue or Darboux-Baire as double
conditions in our results.

While Darboux functions in general may be badly behaved, Darboux functions
of Baire class 1 are more tractable; recall the Kuratowski-Sierpi\'{n}ski
theorem of \S 4. See e.g. [BruL, \S 5], [BruC, \S 6], [CedP], [EvH] for
Darboux functions of Baire class 1, and Marcus [Mar], [GibN1], [GinN2] for
literature and illuminating examples in the study of the Darboux property.%
\newline
\noindent 5. \textit{Qualitative versus quantitative measure theory. }%
Bloom's proof of his theorem used quantitative measure theory. Our proof
replaces this by qualitative measure theory, thus allowing use of
measure-category duality.\newline
\indent The application of Theorem CET above requires the verification of
(wcc), and in the measure case this calls for just enough of the
quantitative aspects to suffice -- see [BinO6, \S 6]. The Baire and measure
cases come together here via the coincidence between measure and metric for
real intervals, cf. \S 5.11.\newline
\noindent 6. \textit{Beyond the reals. }Theorem CET above was conceived to
capture topologically the embedding properties enjoyed by non-negligible
sets under translation as typified by the Steinhaus Theorem, or Sum-Set
Theorem (that $A-A$ has $0$ in its interior). Thus CET\ refers to the
underlying group of homeomorphisms of a space. A more general setting
involves the apparatus of group action on a topological space (see [MilO]).
Here the central result is the Effros Theorem, which may be deduced from
CET-like theorems (see [Ost3]).\newline
\indent The context in our results here is real analysis, as in BGT and
[Blo]. But the natural setting is much more general. One such setting is the
normed groups of [BinO6] (where one has the dichotomy: normed groups are
either topological, or pathological); see [BinO5, Part I] for a development
of slow variation in that context. Other possible settings include
semitopological groups, paratopological groups, etc.; see e.g. [ElfN].

\noindent 7. \textit{Monotone rearrangements. }The theory of monotone
rearrangements is considered in the last chapter of Hardy, Littlewood and P%
\'{o}lya [HarLP] Ch. X. For a function $f$ its distribution function $%
|\{u:f(u)\leq x\}|$ is non-decreasing and so has a non-decreasing inverse
function, called the non-decreasing (briefly, \textit{increasing}) \textit{%
rearrangement} $f_{\uparrow }$ (thus $f$ and $f_{\uparrow }$ have the same
distribution function). Such rearrangements are of great interest, and use,
in a variety of contexts, including\newline
(i) \textit{probability} (Barlow [Bar], Marcus and Rosen [MarR] \S 6.4);%
\newline
(ii) \textit{statistics}: estimation under monotonicity constraints [JanW];%
\newline
(iii) \textit{optimal transport}: in transport problems with $f$ a strategy $%
f_{\uparrow }$ gives the optimal strategy [Vil];\newline
(iv) \textit{analysis}: [Graf, \S 1.4.1] , [HorW], [BerLR].\newline
As we have seen, the uniformity result for the \textit{monotone} case is
quite simple -- much simpler than for the general case -- and also, $\varphi 
$ monotone will typically be clear from context. In the general case we may
aim to replace $\varphi $ by $\varphi _{\uparrow };$ for specific $\varphi $%
, replacing $e$ (in \S 5.2) by $e_{+}$ may well suffice.\newline
\noindent 8. \textit{The class }$\Gamma $\textit{\ }(BGT \S 3.10) consist of
functions \ $f:\mathbb{R}\rightarrow $ $\mathbb{R}$, non-decreasing
right-continuous which, for some measurable $g:\mathbb{R}\rightarrow (0,%
\mathbb{\infty }),$ the \textit{auxiliary function} of \ $f$, 
\[
f(x+ug(x))/f(x)\rightarrow e^{u}\text{\qquad }x\rightarrow \mathbb{\infty }%
\text{\qquad }(\forall u\in \mathbb{R}).
\]%
It turns out that the convergence here is uniform on compact $u$-sets (from $%
f$ being monotone -- as in Th. 1), and hence that $g$ is self-neglecting.%
\newline
The class $\Gamma $ originates in extreme-value theory (EVT) in probability
theory, in connection with de Haan's work on the domain-of-attraction
problem for the Gumbel (double-exponential) extremal law $\Lambda $ ($%
\Lambda (x):=\exp \{-e^{-x}\}).$ See BGT \S 8.13, [BalE].\newline
\noindent 9. \textit{Beurling regular variation. }In a sequel [BinO10] we
explore the consequences of the \textit{Beurling regular variation} property%
\[
f(x+u\varphi (x))/f(x)\rightarrow g(u)\text{\qquad }x\rightarrow \mathbb{%
\infty }\text{\qquad }(\forall u\in \mathbb{R}).\text{\qquad }(BRV)
\]%
We obtain, in particular, a characterization theorem%
\[
g(u)=e^{\rho u}
\]%
for some $\rho .$ We also relax the condition above that $f$ be monotone.%
\newline
\indent Just as Beurling slow variation can be developed in contexts more
generic than the reals (\S 5.6), so too can Beurling regular variation, a
theme that we develop elsewhere [BinO11].\newline
\noindent 10. \textit{Continuous and sequential aspects. }The reader will
have noticed that Beurling slow variation is (like Karamata slow variation)
a continuous-variable property, while the proofs here are by contradiction,
and use sequences (bearing witness to the contradiction). This is a
recurrent theme; see e.g. BGT \S 1.9, [BinO6].\newline
\noindent 11. \textit{The Weil topology and logarithms: lines and
half-lines. }Recall that in connection with Theorem W in \S 1 we mentioned
the additive form with Haar measure $dx$ on the line and the multiplicative
form with Haar measure $dx/x$ on the half-line. The relevant background here
is the Weil topology ([Halm, \S 62], [Wei], cf. [BinO6, Th. 6.10], see also
[HewR]). The relevant metrics involve the logarithmic mapping ($\log
x=\int_{[1,x]}du/u$ and $|[a,b]|_{W}=\log b-\log a)$. This theme underlies
the passage between $\varphi $ and $h=\log \varphi ,$ and is important in
the sequel, [BinO10].\newline
\noindent 12. \textit{Questions. }We close with two questions.\newline
1. Does Bloom's theorem extend to measurable/Baire functions -- that is, can
one omit the Darboux requirement? Does it even extend to Baire-1 functions?%
\newline
2. Are the classes $BSV$, $SN\ $closed under monotone rearrangement?

\begin{center}
\textbf{References}
\end{center}

\noindent \lbrack BalE] A. A. Balkema and P. Embrechts, \textsl{High risk
scenarios and extremes: A geometric approach}. European Math. Soc., Z\"{u}%
rich, 2007. \newline
\noindent \lbrack Bar] M. T. Barlow, Necessary and sufficient conditions for
the continuity of local times of Levy processes. \textsl{Ann. Prob.} \textbf{%
16} (1988), 1389-1427.\newline
\noindent \lbrack BerLR] H. Berestycki and T. Lachand-Robert, Some
properties of monotone rearrangements with applications to elliptic
equations in cylinders. \textsl{Math. Nachrichten} \textbf{266} (2004), 3-19.%
\newline
\noindent \lbrack Bin1] N. H. Bingham, Tauberian theorems and the central
limit theorem. \textsl{Ann. Probab.} \textbf{9} (1981), 221-231. \newline
\noindent \lbrack Bin2] N. H. Bingham, On Euler and Borel summability. 
\textsl{J. London Math. Soc.} (2) \textbf{29} (1984), 141-146. \newline
\noindent \lbrack Bin3] N. H. Bingham, On Valiron and circle convergence. 
\textsl{Math. Z.} \textbf{186} (1984), 273-286. \newline
\noindent \lbrack Bin4] N. H. Bingham, Tauberian theorems for summability
methods of random-walk type. \textsl{J. London Math. Soc.} (2) \textbf{30}
(1984), 281-287.\newline
\noindent \lbrack Bin5] N. H. Bingham, Moving averages. \textsl{Almost
Everywhere Convergence I} (ed. G.A. Edgar \& L. Sucheston) 131-144, Academic
Press, 1989.\newline
\noindent \lbrack Bin6] N. H. Binghams, Szeg\H{o}'s theorem and its
probabilistic descendants. \textsl{Probability Surveys }\textbf{9 }(2012),
287-324.\newline
\noindent \lbrack BinG1] N. H. Bingham and C.M. Goldie, On one-sided
Tauberian conditions. \textsl{Analysis} \textbf{3} (1983), 159-188.\newline
\noindent \lbrack BinG2] N. H. Bingham and C. M. Goldie, Riesz means and
self-neglecting functions. \textsl{Math. Z.} \textbf{199} (1988), 443-454.%
\newline
\noindent \lbrack BinGT] N. H. Bingham, C. M. Goldie and J. L. Teugels, 
\textsl{Regular variation}. 2nd ed., Cambridge University Press, 1989 (1st
ed. 1987). \newline
\noindent \lbrack BinO1] N. H. Bingham and A. J. Ostaszewski, Infinite
combinatorics and the foundations of regular variation. \textsl{J. Math.
Anal. Appl.} \textbf{360} (2009), 518-529. \newline
\noindent \lbrack BinO2] N. H. Bingham and A. J. Ostaszewski, Infinite
combinatorics in function spaces: category methods. \textsl{Publ. Inst.
Math. (Beograd)} (N.S.) \textbf{86 }(100) (2009), 55--73.\newline
\noindent \lbrack BinO3] N. H. Bingham and A. J. Ostaszewski, Beyond
Lebesgue and Baire: generic regular variation. \textsl{Colloq. Mathematicum} 
\textbf{116} (2009), 119-138. \newline
\noindent \lbrack BinO4] N. H. Bingham and A. J. Ostaszewski, Beyond
Lebesgue and Baire II: bitopology and measure-category duality. \textsl{%
Colloq. Math.} \textbf{121} (2010), 225-238.\newline
\noindent \lbrack BinO5] N. H. Bingham and A. J. Ostaszewski, Topological
regular variation. I: Slow variation; II: The fundamental theorems; III:
Regular variation. \textsl{Topology and its Applications} \textbf{157}
(2010), 1999-2013, 2014-2023, 2024-2037. \newline
\noindent \lbrack BinO6] N. H. Bingham and A. J. Ostaszewski, Normed groups:
Dichotomy and duality. \textsl{Dissertationes Math.} \textbf{472} (2010),
138p. \newline
\noindent \lbrack BinO7] N. H. Bingham and A. J. Ostaszewski, Kingman,
category and combinatorics. \textsl{Probability and Mathematical Genetics}
(Sir John Kingman Festschrift, ed. N. H. Bingham and C. M. Goldie), 135-168,
London Math. Soc. Lecture Notes in Mathematics \textbf{378}, CUP, 2010. 
\newline
\noindent \lbrack BinO8] N. H. Bingham and A. J. Ostaszewski, Dichotomy and
infinite combinatorics: the theorems of Steinhaus and Ostrowski. \textsl{%
Math. Proc. Cambridge Phil. Soc.} \textbf{150} (2011), 1-22. \newline
\noindent \lbrack BinO9] N. H. Bingham and A. J. Ostaszewski, Steinhaus
theory and regular variation: De Bruijn and after. \textsl{Indagationes
Mathematicae} (N. G. de Bruijn Memorial Issue), to appear.\newline
\noindent \lbrack BinO10] N. H. Bingham and A. J. Ostaszewski, Uniformity
and self-neglecting functions: II. Beurling regular variation and the class $%
\Gamma ,$ preprint (available at: http://www.maths.lse.ac.uk/Personal/adam/).%
\newline
\noindent \lbrack BinO11] N. H. Bingham and A. J. Ostaszewski, Beurling
regular variation in Banach algebras: asymptotic actions and cocycles, in
preparation.\newline
\noindent \lbrack BinT] N. H. Bingham and G. Tenenbaum, Riesz and Valiron
means and fractional moments. \textsl{Math. Proc. Cambridge Phil. Soc.} 
\textbf{99} (1986), 143-149.\newline
\noindent \lbrack Blo] S. Bloom, A characterization of B-slowly varying
functions. \textsl{Proc. Amer. Math. Soc.} \textbf{54} (1976), 243-250. 
\newline
\noindent \lbrack Boa1] R. P. Boas, \textsl{Entire functions}. Academic
Press, 1954. \newline
\noindent \lbrack Boa2] R. P. Boas, \textsl{A primer of real functions}. 3rd
ed. Carus Math. Monographs 13, Math. Assoc. America, 1981. \newline
\noindent \lbrack Bou] A. Bouziad, The point of continuity property,
neighbourhood assignements and filter convergence. \textsl{Fund. Math.} 
\textbf{218} (2012), 225-242.\newline
\noindent \lbrack BruC] A. M. Bruckner and J. G. Ceder, \textsl{Darboux
continuity}. Jahresbericht d. DMV 67 (1965), 93-117. (Available at:\newline
http://gdz.sub.uni-goettingen.de/dms/load/img/?PPN=\newline
PPN37721857X\_0067\&DMDID=DMDLOG\_0011\&LOGID=LOG\_0011\&\newline
PHYSID=PHYS\_0101)\newline
\noindent \lbrack BruCW] A. M. Bruckner, J. G. Ceder and M. Weiss, \textsl{%
On uniform limits of Darboux functions}. Coll. Math. \textbf{15} (1966),
65-77. \newline
\noindent \lbrack BruL] A. M. Bruckner and J. L. Leonard, Derivatives. 
\textsl{Amer. Math. Monthly}, \textbf{37} (1966), 24-56. \newline
\noindent \lbrack CedP] J. Ceder and T. Pearson, A survey of Darboux Baire 1
functions. \textsl{Real Analysis Exchange }\textbf{9} (1983-84), 179-194.%
\newline
\noindent \lbrack ChaM] K. Chandrasekharan and S. Minakshisundaram, \textsl{%
Typical means}. Oxford University Press, 1952. \newline
\noindent \lbrack ElfN] A. S. Elfard and P. Nickolas, On the topology of
free paratopological groups. \textsl{Bull. London Math. Soc. }\textbf{44}
(2012), 1103-1115.\newline
\noindent \lbrack EvH] M. J. Evans and P. D. Humke, Revisiting a century-old
characterization of Baire one Darboux functions. \textsl{Amer. Math. Monthly 
}\textbf{116} (2009), 451-455.\newline
\noindent \lbrack FeiS] H. G. Feichtinger and H. J. Schmeisser, Weighted
versions of Beurling's slowly varying functions. \textsl{Math. Ann.} \textbf{%
275} (1986), 353-363. \newline
\noindent \lbrack GibN1] R. G. Gibson and T. Natkaniec, Darboux like
functions. \textsl{Real Analysis Exchange} \textbf{22.2} (1996/97), 492-453. 
\newline
\noindent \lbrack GibN2] R. G. Gibson and T. Natkaniec, Darboux-like
functions. Old problems and new results. \textsl{Real Analysis Exchange} 
\textbf{24.2} (1998/99), 487-496. \newline
\noindent \lbrack Graf] L. Grafakos, \textsl{Classical Fourier Analysis},
2nd ed., Graduate Texts in Mathematics 249, Springer, 2008.\newline
\noindent \lbrack Halm] P. R. Halmos. \textsl{Measure theory, }(1955) Van
Nostrand, 1950 (Grad. Texts in Math. 18, Springer, 1970).\newline
\noindent \lbrack Halp1] I. Halperin, On the Darboux property. \textsl{%
Pacific J. Math.} \textbf{5} (1955), 703-705. \newline
\noindent \lbrack Halp2] I. Halperin, Discontinuous functions with the
Darboux property. \textsl{Canad. Math. Bull.} \textbf{2} (1959), 111-118. 
\newline
\noindent \lbrack Har] G. H. Hardy, \textsl{Divergent series}. Oxford
University Press, 1949. \newline
\noindent \lbrack HarLP] G. H. Hardy, J. E. Littlewood and G. P\'{o}lya. 
\textsl{Inequalities}. 2nd ed., CUP, 1952 (1st. ed. 1934)\newline
\noindent \lbrack HorW] A. Horsley and A. J. Wr\'{o}bel, The Mackey
continuity of the monotone rearrangement. \textsl{Proc. Amer. Math. Soc.} 
\textbf{97} (1986), 626-628.\newline
\noindent \lbrack HewR] E. Hewitt and K. A. Ross, \textsl{Abstract Harmonci
Analysis, I Structure of topological groups, integration theory, group
representations}$.$ Grundl. math. Wiss. 115, Springer, 1963.\newline
\noindent \lbrack JanW] H. K. Jankowski and J. A. Wellner, Estimation of
discrete monotone distributions. \textsl{Electronic J. Stat.} \textbf{3}
(2009), 1567-1605.\newline
\noindent \lbrack Kar] J. Karamata, Sur un mode de croissance r\'{e}guli\`{e}%
re des fonctions. \textsl{Mathematica (Cluj)} \textbf{4} (1930), 38-53.%
\newline
\noindent \lbrack Kech] A. S. Kechris, \textsl{Classical Descriptive Set
Theory.} Grad. Texts in Math. 156, Springer, 1995.\newline
\noindent \lbrack Koo] P. Koosis, \textsl{The logarithmic integral}. I 2nd
ed. CUP, 1998 (1s ed. 1988), II, CUP, 1992. \newline
\noindent \lbrack Kor] J. Korevaar, \textsl{Tauberian theorems: A century of
development}. Grundl. math. Wiss. \textbf{329}, Springer, 2004. \newline
\noindent \lbrack KorvAEdB] J. Korevaar, T. van Aardenne-Ehrenfest and N. G.
de Bruijn: A note on slowly oscillating functions. \textsl{Nieuw Arch.
Wiskunde} \textbf{23} (1949), 77-86. \newline
\noindent \lbrack Kur] C. Kuratowski, \textsl{Topologie}. Monografie Mat. 20
(4th. ed.), PWN Warszawa 1958 [K. Kuratowski, \textsl{Topology}. Translated
by J. Jaworowski, Academic Press-PWN 1966]. \newline
\noindent \lbrack KurS] K. Kuratowski and W. Sierpi\'{n}ski, Sur les
fonctions de classe I et les ensembles connexes ponctiformes. \textsl{Fund.
Math.} \textbf{3}, (1922) 303-313. \newline
\noindent \lbrack LeeTZ] P.-Y. Lee, W.-K. Tang and D. Zhao, An equivalent
definition of functions of first Baire class. \textsl{Proc. Amer. Math. Soc.}
\textbf{129} (2001), 2273-2275.\newline
\noindent \lbrack Mar] S. Marcus, Functions with the Darboux property and
functions with connected graphs. \textsl{Math. Annalen.} \textbf{141}
(1960), 311-317. \newline
\noindent \lbrack MarR] M. B. Marcus and J. Rosen, \textsl{Markov processes,
Gaussian processes and local time}. Cambridge Studies in Adv. Math., 100,
CUP, 2006.\newline
\noindent \lbrack Mat] W. Matuszewska, On a generalisation of regularly
increasing functions. \textsl{Studia Math.} \textbf{24} (1962), 271-276. 
\newline
\noindent \lbrack MilO] H. I. Miller and A. J. Ostaszewski, Group actions
and shift-compactness. \textsl{J. Math. Anal. Appl.} \textbf{392} (2012),
23-39.\newline
\noindent \lbrack Moh] T. T. Moh, On a general Tauberian theorem. \textsl{%
Proc. Amer. Math. Soc.} \textbf{36} (1972), 167-172.\newline
\noindent \lbrack Nat] I. P. Natanson,\textsl{\ Theory of functions of a
real variable. }Vol. I,II, Frederick Ungar, 1960. \newline
\noindent \lbrack Ost1] A. J. Ostaszewski, Analytic Baire spaces. \textsl{%
Fund. Math.} \textbf{217} (2012), 189-210. \newline
\noindent \lbrack Ost2] A. J. Ostaszewski, Almost completeness and the
Effros Theorem in normed groups. \textsl{Topology Proceedings} \textbf{41}
(2013), 99-110.\newline
\noindent \lbrack Ost3] A. J. Ostaszewski, Shift-compactness in almost
analytic submetrizable Baire groups and spaces, invited survey article. 
\textsl{Topology Proceedings} \textbf{41} (2013), 123-151.\newline
\noindent \lbrack Ost4] A. J. Ostaszewski, On the Effros Open Mapping
Theorem -- separable and non-separable, preprint (available at:
www.maths.lse.ac.uk/Personal/adam).\newline
\noindent \lbrack Oxt] J. C. Oxtoby, \textsl{Measure and category.} 2nd ed.,
Grad. Texts Math. \textbf{2}, Springer, 1980. \newline
\noindent \lbrack Pet] G. E. Petersen, Tauberian theorems for integrals II. 
\textsl{J. London Math. Soc.} \textbf{5} (1972), 182-190. \newline
\noindent \lbrack PolS] G. P\'{o}lya and G. Szeg\H{o}, \textsl{Aufgaben und
Lehrs\"{a}tze aus der Analysis} Vol. I. Grundl. math. Wiss. XIX, Springer,
1925. \newline
\noindent \lbrack PorWBW] W. Poreda, E. Wagner-Bojakowska and W. Wilczy\'{n}%
ski, A category analogue of the density topology. \textsl{Fund. Math.} 
\textbf{125} (1985), 167-173. \newline
\noindent \lbrack Rud] W. Rudin, \textsl{Functional analysis}. 2nd ed.
McGraw-Hill, 1991 (1st ed. 1973). \newline
\noindent \lbrack Sol] S. Solecki, Decomposing Borel sets and functions and
the structure of Baire class 1 functions. \textsl{J. Amer. Math. Soc. }%
\textbf{11.3} (1998), 521-550.\newline
\noindent \lbrack Ten] G. Tenenbaum, Sur le proc\'{e}d\'{e} de sommation de
Borel et la r\'{e}partition des facteurs premiers des entiers. \textsl{%
Enseignement Math.} \textbf{26} (1980), 225-245. \newline
\noindent \lbrack Vil] C. Villani, \textsl{Topics in optimal transportation}%
. Grad. Studies in Math. 58, Amer. Math. Soc., 2003.\newline
\noindent \lbrack Wei] A. Weil, \textsl{L'int\'{e}gration dans les groupes
topologiques et ses applications}. Actual. Sci. Ind. 689, Hermann, Paris,
1940 (republished, Princeton Univ. Press, 1941).\textsl{\newline
}\noindent \lbrack Wid] D. V. Widder, \textsl{The Laplace Transform}.
Princeton, 1941.\newline
\noindent \lbrack Wie] N. Wiener, Tauberian theorems. \textsl{Acta Math.} 
\textbf{33} (1932), 1-100 (reprinted in N. Wiener, \textsl{Generalized
harmonic analysis and Tauberian theorems}. MIT Press, 1964). \newline

\bigskip

\noindent Mathematics Department, Imperial College, London SW7 2AZ;
n.bingham@ic.ac.uk \newline
Mathematics Department, London School of Economics, Houghton Street, London
WC2A 2AE; A.J.Ostaszewski@lse.ac.uk 

\end{document}